\documentclass{article}

\usepackage{amsmath,amsfonts,amsthm,amssymb,graphicx}
\usepackage[all,2cell,ps]{xy}

\theoremstyle{plain}
\newtheorem{thm}{Theorem}[section]
\newtheorem{lem}[thm]{Lemma}
\newtheorem{prop}[thm]{Proposition}
\newtheorem{cor}[thm]{Corollary}

\theoremstyle{definition}

\newtheorem*{rem}{Remark}

\theoremstyle{remark}

\newcommand{\bs}{\backslash}

\DeclareMathOperator{\U}{U}

\newcommand{\Z}{\mathbb Z}

\newcommand{\R}{\mathbb R}

\newcommand{\Q}{\mathbb Q}

\newcommand{\gam}{\gamma}

\newcommand{\Lam}{\Lambda}
\newcommand{\A}{\mathbb{A}}
\newcommand{\C}{\mathbb C}

\newcommand{\Gam}{\Gamma}

\newcommand{\conj}{\overline}

\newcommand{\mc}{\mathcal}
\newcommand{\del}{\delta}

\DeclareMathOperator{\SL}{SL}
\DeclareMathOperator{\PSL}{PSL}
\DeclareMathOperator{\GL}{GL}

\DeclareMathOperator{\SO}{SO}

\DeclareMathOperator{\SU}{SU}

\DeclareMathOperator{\Sp}{Sp}

\DeclareMathOperator{\Spin}{Spin}

\newcommand{\mf}{\mathfrak}

\newenvironment{pf}{\begin{proof}}{\end{proof}}

\title{On the number of ends of rank one locally symmetric spaces}
\author{Matthew Stover\footnote{Partially supported by NSF RTG grant DMS 0602191.} \\ \small{University of Michigan}\\ \small{\textsf{stoverm@umich.edu}}}
\date{\today}

\begin{document}

\maketitle

\begin{abstract}
Let $Y$ be a noncompact rank one locally symmetric space of finite volume. Then $Y$ has a finite number $e(Y) > 0$ of topological ends. In this paper, we show that for any $n \in \mathbb N$, the $Y$ with $e(Y) \leq n$ that are arithmetic fall into finitely many commensurability classes. In particular, there is a constant $c_n$ such that $n$-cusped arithmetic orbifolds do not exist in dimension greater than $c_n$. We make this explicit for one-cusped arithmetic hyperbolic $n$-orbifolds and prove that none exist for $n \geq 30$.
\end{abstract}


\section{Introduction}\label{sec:intro}


Let $X$ be real, complex, quaternionic hyperbolic space, or the Cayley hyperbolic plane and $G$ its orientation-preserving isometry group. Let $\Gam < G$ be a lattice and $\Gam \bs X$ be the associated locally symmetric space. Throughout this paper we assume that $\Gam \bs X$ is noncompact, i.e., that $\Gam$ is a nonuniform lattice in $G$. Then $\Gam \bs X$ has a finite number of topological ends $e(\Gam \bs X) > 0$. The purpose of this paper is to prove the following theorem.


\begin{thm}\label{thm:intro}
Fix $n > 0$. There are only finitely many commensurability classes of arithmetic rank one locally symmetric spaces $Y$ with $e(Y) = n$.
\end{thm}


Some finiteness results were previously known for a fixed rank one symmetric space. For the hyperbolic plane, the modular group $\PSL_2(\Z)$ determines the unique commensurability class of nonuniform arithmetic lattices in $\PSL_2(\R)$. For hyperbolic $3$-space see Chinburg, Long, and Reid \cite{Chinburg--Long--Reid}, and for the complex hyperbolic plane see \cite{Stover}. Every lattice in $\mathrm{PSp}(n, 1)$ or $\mathrm F_4^{(-20)}$ is arithmetic \cite{Corlette, Gromov--Schoen}, so the arithmetic assumption is superfluous and we have the following.


\begin{cor}\label{cor:quaternion}
For every $n > 0$ there are, up to commensurability, only finitely many $n$-cusped quaternionic and Cayley hyperbolic orbifolds of finite volume.
\end{cor}


In the language of algebraic groups, the arithmetic groups that we consider are those which have $\Q$- and $\R$-rank one. While the $\R$-rank one assumption is natural from a geometric perspective, since this corresponds to the negatively curved symmetric spaces, it is also natural from the point of view of algebraic groups. For example, all finite volume locally symmetric spaces of $\Q$-rank at least $2$ (e.g., $\SL_n(\Z)$ for $n \geq 3$) are one-ended, and many remain `$1$-cusped' under various alternative definitions of a cusp, like transitivity of the action on $\Q$-isotropic flags. See \cite{Borel--Ji} for various interpretations in higher rank. For $\Q$-rank $1$ lattices in semisimple groups of higher $\R$-rank, one encounters lattices like $\SL_2(\mathcal O_F)$, where $\mathcal O_F$ is the ring of integers in an arbitary algebraic number field. The number of ends of the locally symmetric space associated with $\SL_2(\mathcal O_F)$ equals the so-called class number of $F$, and whether or not there are infinitely many number fields of bounded class number is widely expected to be true, but remains one of the outstanding open problems in number theory.

We now introduce notation that we will use throughout the paper. Since $G$ has real rank one, there is a unique conjugacy class of parabolic subgroups. Let $P$ denote one such subgroup. Choose a maximal $\R$-split torus $S \subset P$, and let $Z$ be the centralizer of $S$ in $G$. Then there is a unipotent subgroup $U \subset P$ so that $P$ is the semidirect product of $U$ with $Z$.

Since $P$ is the stabilizer in $G$ of a point on the ideal boundary $X_\infty$ of $X$ and $G$ acts transitively on the boundary, $X_\infty$ is naturally identified with the coset space $G / P$. The ends of $\Gam \bs X$ are in one-to-one correspondence with the $\Gam$-conjugacy classes of parabolic subgroups of $\Gam$. In other words, ends correspond to $\Gam$-orbits of those $g P \in G / P$ such that $\Gam \cap g P g^{-1}$ is a cocompact lattice in $g P g^{-1}$. This leads to the following interpretation of the ends of $\Gam \bs X$ when $\Gam$ is arithmetic.

Let $\Gam < G$ be a nonuniform arithmetic lattice. Then there is an absolutely almost simple simply connected $\Q$-algebraic group $\mc G$ such that the lift of $\Gam$ to $\mc G(\R)$ under a central isogeny $\mc G(\R) \to G$ is commensurable with the group of integral points $\mc G(\Z)$ determined by a representation of $\mc G$ into $\GL_N(\Q)$. See \cite{Borel1, Borel3}. We describe the construction of these lattices in $\S$\ref{sec:arithmetic}.

Since $\mc G$ has $\Q$-rank one, there is a unique conjugacy class of $\Q$-parabolic subgroups. Let $\mc P$ be one, and choose a maximal $\Q$-split torus $\mc S \subset \mc P$. Then $\Gam$ acts on the complete variety $(\mc G / \mc P)_\Q = \mc G(\Q) / \mc P(\Q)$. Since $\mc P(\R)$ contains the center of $\mc G(\R)$, we can identify $\mc G(\R) / \mc P(\R)$ with $X_\infty$. Parabolic subgroups of $\Gam$ are commensurable with a lattice in $\mc P(\Q)$, so
\begin{equation}\label{eq:cuspcount}
e(\Gam \bs X) = \# \big( \Gam \bs \mc G(\Q) / \mc P(\Q) \big).
\end{equation}
The focus of this paper is on the right-hand side of \eqref{eq:cuspcount}.

We begin in $\S$\ref{sec:principal} by studying the ends of $\Gam \bs X$ when $X$ is the image in $G$ of a so-called \emph{principal arithmetic} subgroup $\Gam_{K_f}$ of $\mc G(\Q)$ defined by a coherent compact open subgroup $K_f$ of $G(\A_f)$, where $\A_f$ denotes the finite adeles of $\Q$. When $K_f$ is \emph{special} (see \cite{Tits2}), then we can give an exact formula for $e(\Gam_{K_f} \bs X)$ using results of Borel \cite{Borel2}. In other cases, we only obtain a lower bound. See Proposition \ref{prop:cusp-centralizer}.

To prove Theorem \ref{thm:intro}, it suffices to consider maximal arithmetic lattices in $G$. By Proposition 1.4 in \cite{Borel--Prasad}, every maximal arithmetic lattice is the normalizer in $G$ of some principal arithmetic lattice. This is analyzed in $\S$\ref{sec:maximal}, where we complete the proof of Theorem \ref{thm:intro}. In $\S$\ref{sec:onecusp}, we apply our techniques to give an explicit bound for one-cusped arithmetic hyperbolic $n$-orbifolds.


\begin{thm}\label{thm:intro1cusp}
One-cusped arithmetic hyperbolic $n$-orbifolds do not exist for any $n \geq 30$.
\end{thm}


It is known that there are hyperbolic reflection groups that determine one-cusped arithmetic hyperbolic $n$-orbifolds for all $n \leq 9$ \cite{Hild}. We close the paper by constructing one-cusped hyperbolic $n$-orbifolds for $n = 10, 11$. There may be examples for $12 \leq n \leq 29$ related to definite rational quadratic forms with few classes in their spinor genus; such quadratic forms do not seem to be classified, so we do not know if Theorem \ref{thm:intro1cusp} is sharp.


\subsubsection*{Acknowledgments}

I thank Gopal Prasad for some helpful suggestions while I was completing this paper.


\section{Arithmetic subgroups of rank one groups}\label{sec:arithmetic}


In this section, we describe the nonuniform arithmetic lattices in simply connected Lie groups of $\R$-rank one. See \cite{Tits1} for the full classification. This naturally breaks up into three cases: hyperbolic space, complex and quaternionic hyperbolic space, and the Cayley hyperbolic plane.


\subsection{Hyperbolic space}\label{subsec:real}


The simply connected form of the isometry group of hyperbolic $n$-space is the group $\Spin(n, 1)$ \cite[$\S$27.4B]{Knus--etal}, which is the double-cover of $\SO(n, 1)$. Note that we have exceptional isomorphisms
\[
\Spin(2, 1) \cong \SL_2(\R),
\]
\[
\Spin(3, 1) \cong \SL_2(\C)
\]
with the more familiar groups acting on hyperbolic $2$- and $3$-space.

All nonuniform arithmetic lattices in $\Spin(n, 1)$ are determined as follows. Let $q$ be an isotropic nondegenerate quadratic form on $\Q^{n + 1}$ of signature $(n, 1)$ and $\mc G = \Spin(q)$. Recall that a quadratic form is isotropic if there is a nonzero vector $v \in \Q^{n + 1}$ so that $q(v) = 0$. Then $\mc G(\R) \cong \Spin(n, 1)$, and every nonuniform arithmetic lattice in $\Spin(n, 1)$ is commensurable with $\mc G(\Z)$ for some $\mc G$ as above. We note that there are constructions of arithmetic lattices for all odd $n$ that do not use quadratic forms, but these constructions do not lead to nonuniform lattices. See \cite{Vinberg--Shvartsman} for further details.

We now describe the $\Q$-split tori of $\mc G$ and their centralizers, since they are crucial throughout this paper. A maximal $\Q$-split torus $\mc S$ of $\mc G$ is isomorphic to the multiplicative group $\mathbb G_m$ over $\Q$. Since $q$ is isotropic, we can find a basis for $\Q^{n + 1}$ such that $q$ has matrix
\[
Q = \begin{pmatrix}
0 & 0 & 1 \\
0 & Q' & 0 \\
1 & 0 & 0
\end{pmatrix},
\]
where $Q'$ is the matrix of an anisotropic (i.e., not isotropic) quadratic form $q'$ on $\Q^{n - 1}$. That is, $q = q_0 \oplus q'$, where $q_0$ is a hyperbolic plane. Under this basis, the image $\conj{\mc S}$ of $\mc S$ in the special orthogonal group $\SO(q)$ is the set of matrices of the form
\begin{equation}\label{eq:torus-form}
\begin{pmatrix}
x & 0 & 0 \\
0 & \mathrm I_{n - 1} & 0 \\
0 & 0 & x^{-1}
\end{pmatrix},
\end{equation}
where $x \in \Q^*$ and $\mathrm I_{n - 1}$ is the $(n - 1) \times (n - 1)$ identity matrix.

We also need to understand the centralizer $Z(\mc S)$ of $\mc S$ in $\mc G$. The centralizer $Z(\conj{\mc S})$ of $\conj{\mc S}$ in $\SO(q)$ is the group of elements
\[
\begin{pmatrix}
x & 0 & 0 \\
0 & A & 0 \\
0 & 0 & x^{-1}
\end{pmatrix}
\]
such that $x \in \Q^*$ and $A \in \SO(q')$. Then $Z(\mc S)$ is the lift of this group to $\Spin(q)$, and the quotient of $Z(\mc S)$ by $\mc S$ is the group $\mathrm{Pin}(q')$, which contains $\Spin(q')$ as an index $2$ subgroup.


\subsection{Complex and quaternionic hyperbolic space}\label{subsec:C-H}


The simply connected forms of the isometry group of complex and quaternionic hyperbolic $n$-spaces are $\SU(n, 1)$ and $\Sp(n, 1)$, respectively \cite[$\S$24.7B]{Knus--etal}. For complex hyperbolic space, let $D$ be an imaginary quadratic extension of $\Q$. For quaternionic hyperbolic space, let $D$ be a definite quaternion algebra with center $\Q$, i.e., one such that $D \otimes_\Q \R$ is isomorphic to Hamilton's quaternions $\mathbb H$. We then have an involution $\tau : D \to D$ given by the nontrivial Galois automorphism when $D$ is an imaginary quadratic field and quaternion conjugation when $D$ is a quaternion algebra.

Let $h$ be an isotropic nondegenerate $\tau$-hermitian form on $D^{n + 1}$. Then $h$ is a $\tau$-symmetric matrix in $\GL_{n + 1}(\C)$ or $\GL_{n + 1}(\mathbb H)$, so it has real eigenvalues. Therefore, the signature of $h$ makes sense, and we assume that $h$ has signature $(n, 1)$. If $\mc G$ is the special automorphism group of $h$, then $\mc G(\R)$ is isomorphic to $\SU(n, 1)$ when $D$ is imaginary quadratic and $\Sp(n, 1)$ when $D$ is a quaternion algebra. Any nonuniform arithmetic lattice in $\SU(n, 1)$ or $\Sp(n, 1)$ is commensurable with $\mc G(\Z)$ for some $\mc G$ as above. Again, there are other constructions of cocompact lattices, but these algebraic groups suffice for the nonuniform lattices.

We again describe some facts about centralizers of $\Q$-split tori that we will need later. As in the hyperbolic case, the maximal $\Q$-split torus $\mc S$ is isomorphic to the multiplicative group of $\Q$, and we can choose a basis for $D^{n + 1}$ for which
\[
h = \begin{pmatrix}
0 & 0 & 1 \\
0 & h' & 0 \\
1 & 0 & 0
\end{pmatrix},
\]
where $h'$ is an anisotropic $\tau$-hermitian form on $D^{n - 1}$. Then $\mc S$ is realized as matrices exactly the same as \eqref{eq:torus-form} and the centralizer of $\mc S$ now consists of those matrices
\[
(x, A) = \begin{pmatrix}
x & 0 & 0 \\
0 & A & 0 \\
0 & 0 & \tau(x)^{-1}
\end{pmatrix}
\]
such that $x \in D^*$, $A$ is in the unitary group of $h'$ (not the special unitary group), and $x \tau(x)^{-1} \det(A) = 1$.

We claim that $Z(\mc S) / \mc S$ is isomorphic to $\U(h')$ (as $\Q$-algebraic groups) under the map $(x, A) \mapsto A$. The kernel of this map is clearly $\mc S$, so it suffices to show that this map is onto. That is, given $A \in \U(h')$, we must show that there exists $x \in D^*$ such that $(x, A) \in \SU(h)$. That is, we need to know that there exists $x \in D^*$ such that $x^{-1} \tau(x) = \det(A)$. This follows immediately from Hilbert's Theorem 90, which holds for both imaginary quadratic fields and quaternion algebras \cite[$\S$29.A]{Knus--etal}.


\subsection{The Cayley hyperbolic plane}\label{subsec:O}


See \cite{Allcock} for a more detailed description of lattices in $\mathrm F_4^{(-20)}$. Let $C$ be a Cayley algebra over $\Q$ with involution $\tau$ and $h$ be a $\tau$-symmetric element of $\GL_3(C)$. The automorphisms of $h$ with reduced norm $1$ form an algebraic group $\mc G$ that is simply connected with $\mc G(\R) \cong \mathrm F_4^{(-20)}$. One can also realize this as the automorphisms of an exceptional Jordan algebra. The $\Q$-split torus of $\mc G$ again has the form \eqref{eq:torus-form}, and $Z(\mc S) / \mc S$ is isomorphic over $\Q$ to the group of elements in $C$ with reduced norm $1$.


\section{Principal arithmetic lattices}\label{sec:principal}


We begin with some general results. Let $\A_k$ be the adeles of the number field $k$ and $\A_{k, f}$ the finite adeles. We suppress the $k$ when $k = \Q$.

See \cite{Borel1} for the basic theory of algebraic groups over number fields. Let $\mc G$ be an absolutely almost simple and simply connected $k$-algebraic group and $\mc H$ be a $k$-parabolic subgroup. These assumptions ensure that $\mc G$ has strong approximation, i.e., that $\mc G(k)$ is dense in $\mc G(\A_{k,f})$ \cite[Thm.~7.12]{Platonov--Rapinchuk}. If $K_f \subset \mc G(\A_{k, f})$ is an open compact subgroup, set $K_f^\infty = \mc G(k \otimes \R) \times K_f \subset \mc G(\A)$. Then
\begin{equation}\label{eq:strong-approx}
K_f^\infty \mc G(k) = \mc G(\A).
\end{equation}

Let $K_f$ be an open compact subgroup of $\mc G(\A_{k,f})$ and set $L_f = \mc H(\A_{k,f}) \cap K_f$. Then $\Gam_{K_f} = \mc G(k) \cap K_f$ is a lattice in $\mc G(k \otimes \R)$, and we are interested in the quantity
\begin{equation}\label{eq:h-cusps}
e_{\mc H}(\Gam_{K_f}) = \# \left(\Gam_{K_f} \bs \mc G(k) / \mc H(k) \right).
\end{equation}
When $k = \Q$ and $\mc G(\R)$ has real rank one, then $e_{\mc H}(\Gam_{K_f}) = e(\Gam_{K_f} \bs X)$, where $X$ is the symmetric space associated with $\mc G(\R)$.

In \cite[Prop.~7.5]{Borel2}, Borel relates $e_{\mc H}(\Gam_{K_f})$ to the so-called \emph{class number} of $\mc H$ with respect to $L_f$, which is the number
\begin{equation}\label{eq:h-class}
c(\mc H, L_f) = \# \left( L_f^\infty \bs \mc H(\mathbb A_k) / \mc H(k) \right).
\end{equation}
Also see Chapters 5 and 8 of \cite{Platonov--Rapinchuk}. Since we are restating Borel's results in different language, and because one direction of his proof works greater generality than his stated assumptions, we give a complete proof of \cite[Prop.~7.5]{Borel2} in the next two lemmas. The first step is the following general fact.


\begin{lem}\label{lem:cusp-bound}
Let $\mc G$ be an algebraic group over the number field $k$ and $\mc H$ a $k$-parabolic subgroup. Suppose that $K_f$ is an open compact subgroup of $\mc G(\A_{k, f})$ such that $K_f^\infty \mc G(k) = \mc G(\A_k)$. Let $L_f = \mc H(\A_{k, f}) \cap K_f$ and $\Gam_{K_f} = G(k) \cap K_f$. Then
\begin{equation}\label{eq:cusp-bound}
e_{\mc H}(\Gam_{K_f}) \geq c(\mc H, L_f).
\end{equation}
\end{lem}


\begin{pf}
Given $h_1, h_2 \in \mc H(\A_k)$, there exist $k_1, k_2 \in K_f^\infty$ and $g_1, g_2 \in \mc G(k)$ such that $h_j = k_j g_j$, $j = 1, 2$. Now, suppose that there exists $\gam \in \Gam_{K_f}$ and $h \in \mc H(k)$ such that $g_1 = \gam g_2 h$. Then
\[
h_1 = k_1 g_1 = k_1 \gam g_2 h = (k_1 \gam k_2^{-1}) h_2 h \in K_f^\infty \mc H(\A) \mc H(k).
\]
Since $h_1 (h_2 h)^{-1} \in \mc H(\A_k)$ and $k_1 \gam k_2^{-1} \in K_f^\infty$, we must have $k_1 \gam k_2^{-1} \in L_f^\infty$. Therefore, if $g_1$ and $g_2$ are in the same $\Gam_{K_f}, \mc H(k)$ double coset of $\mc G(k)$, then $h_1$ and $h_2$ are in the same $L_f^\infty, \mc H(k)$ double coset of $\mc H(\A_k)$.

Conversely, suppose $h_1 = \ell h_2 h$ for some $\ell \in L_f^\infty$ and $h \in \mc H(k)$. Then
\[
g_1 = k_1^{-1} h_1 = (k_1^{-1} \ell k_2) g_2 h \in K_f^\infty \mc G(k) \mc H(k).
\]
Since $g_1, g_2, h \in \mc G(k)$, it follows that $k_1^{-1} \gam k_2 \in \Gam_{K_f}$. Therefore, if $h_1$ and $h_2$ are in the same $L_f^\infty, \mc H(k)$ double coset of $\mc H(\A_k)$, then $g_1$ and $g_2$ are in the same $\Gam_{K_f}, \mc H(k)$ double coset of $\mc G(k)$.

It follows that there is a well-defined and injective set map from the finite set $L_f^\infty \bs \mc H(\A_k) / \mc H(k)$ into the finite set $\Gam_{K_f} \bs \mc G(k) / \mc H(k)$. This proves the lemma.
\end{pf}


Note that we did not use one of Borel's assumptions: `$\mathrm G_\mf p = \mathrm G_{\mf o_\mf p}.\mathrm H_\mf p$ for every $\mf p \in \mathrm P$'. In our language, this assumption becomes $\mc G(\A_{k,f}) = K_f \mc H(\A_{k,f})$. When this holds, we say that $\mc G$ has an \emph{Iwasawa decomposition} with respect to $K_f$ and $\mc H$. For example, $\mc G$ has an Iwasawa decomposition when $K_f$ is a coherent product of parahoric subgroups and the $v$-adic component of $K_f$ is maximal and special for every nonarchimedean place $v$ of $k$ \cite[$\S$3.3.2]{Tits2}. The following completes our proof of \cite[Prop.~7.5]{Borel2}.


\begin{lem}\label{lem:cusp-iwa}
With the same assumptions and notation as Lemma \ref{lem:cusp-bound}, suppose that $\mc G$ has an Iwasawa decomposition with respect to $K_f$ and $\mc H$. Then $e_{\mc H}(\Gam_{K_f}) = c(\mc H, L_f)$.
\end{lem}


\begin{pf}
For any $g_1, g_2 \in \mc G(k)$, choose $k_1, k_2 \in K_f^\infty$ and $h_1, h_2 \in \mc H(\mathbb A)$ so that $g_j = k_j h_j$, $j = 1, 2$, under the Iwasawa decomposition of $\mc G$ with respect to $K_f$ and $\mc H$. Note that, by definition of $K_f^\infty$, extending this from $\A_{k, f}$ to $\A_k$ is trivial.

If $g_1 = \gam g_2 h$ for some $\gam \in \Gam_{K_f}$ and $h \in \mc H(k)$, then
\[
h_1 = (k_1^{-1} \gam k_2) h_2 h \in K_f^\infty \mc H(\A_k) \mc H(k).
\]
It follows that $(k_1^{-1} \gam k_2) \in L_f^\infty$ and so $h_1$ and $h_2$ have the same image in $L_f^\infty \bs \mc H(\A_k) / \mc H(k)$. If $h_1$ and $h_2$ lie in the same $L_f^\infty, \mc H(k)$ double coset of $\mc H(\A_k)$, a similar computation shows that $g_1$ and $g_2$ have the same image in $\Gam_{K_f} \bs \mc G(k) / \mc H(k)$. This proves that
\[
e_{\mc H}(\Gam_{K_f}) \leq c(\mc H, L_f),
\]
so the two are equal by Lemma \ref{lem:cusp-bound}.
\end{pf}


We also need the following analogue of \cite[Prop.~2.4]{Borel2}.


\begin{lem}\label{lem:central-ext}
Suppose that
\[
1 \to \mc C \to \mc G \overset{\pi}{\to} \mc H \to 1
\]
is a central exact sequence of $k$-algebraic groups. Let $K_f \subset \mc G(\A_{k, f})$ be an open compact subgroup and $L_f = \pi(K_f)$. Then $c(\mc G, K_f) \geq c(\mc H, L_f)$. Moreover, if $c(\mc C, \mc C(\A_{k, f}) \cap K_f) = 1$ then $c(\mc G, K_f) = c(\mc H, L_f)$.
\end{lem}


\begin{pf}
By assumption, we have a natural surjective map
\[
\hat \pi : K_f^\infty \bs \mc G(\A_k) / \mc G(k) \to L_f^\infty \bs \mc H(\A_k) / \mc H(k).
\]
The first statement follows immediately.

We must show that $\hat \pi$ is injective when $c(\mc C, \mc C(\A_{k, f}) \cap K_f) = 1$. Suppose that $g_1, g_2 \in \mc G(\A_k)$ and $\pi(g_1) = x \pi(g_2) y$ for some $x \in L_f^\infty$ and $y \in \mc H(k)$. Then we have $\tilde x \in K_f^\infty$ and $\tilde y \in \mc G(k)$ such that $g_1^{-1} \tilde x g_2 \tilde y \in \mc C(\A_k)$. However,
\[
\mc C(\A_k) = \left( \mc C(\A_k) \cap K_f^\infty \right) \mc C(k),
\]
so there exist $x_1 \in \mc C(\A_k) \cap K_f^\infty$ and $y_1 \in \mc C(k)$ such that
\[
g_1^{-1} \tilde x g_2 \tilde y = x_1 y_1.
\]
Since $\mc C$ is central, we get
\[
g_1 = (x_1^{-1} x) g_2 (y y_1^{-1}) \in K_f^\infty g_2 \mc G(k).
\]
Thus $\hat \pi$ is injective.
\end{pf}


Now, we return to the case where $k = \Q$ and $\mc G(\R)$ is rank one. Let $\mc P = \mc H$ be a $\Q$-parabolic subgroup of $\mc G$. It is unique up to conjugacy \cite[Prop.~21.12]{Borel1}. Let $\mc S \subset \mc P$ be a maximal $\Q$-split torus. Then $\mc P$ is a semidirect product $Z(\mc S) \mc U$, where $Z(\mc S)$ is the centralizer of $\mc S$ in $\mc G$ and $\mc U$ is unipotent.

We call an open compact subgroup $K_f \subset \mc G(\A_f)$ \emph{coherent} if it is defined by a coherent collection of parahoric subgroups of $\mc G(\Q_p)$ for all $p$; see \cite{Borel--Prasad} or Appendix A of \cite{Platonov--Rapinchuk}. Let $G$ be the orientation preserving isometry group of the symmetric space $X$ and $\Gam < G$ be a nonuniform lattice. We say that $\Gam$ is a \emph{principal arithmetic lattice} if there is an absolutely almost simple and simply connected $\Q$-group $\mc G$ and a coherent open compact subgroup $K_f \subset \mc G(\A_f)$ such that $\Gam$ is the image in $G$ of $\Gam_{K_f} = K_f \cap \mc G(\Q)$ under a central isogeny $\mc G(\R) \to G$. The following allows us to further refine the conclusions of Lemma \ref{lem:cusp-iwa} for principal arithmetic lattices.


\begin{prop}\label{prop:cusp-centralizer}
Let $\mc G$ be a $\Q$ algebraic group of real and $\Q$-rank one and $\mc P$ be a $\Q$-parabolic subgroup with $\mc S \subset \mc P$ a maximal $\Q$-split torus and $Z(\mc S)$ the centralizer of $\mc S$ in $\mc G$. Let $K_f \subset \mc G(\A_f)$ be an open compact subgroup determined by a coherent product of parahoric subgroups. Set:
\[
L_f = K_f \cap \mc H(\A_f),
\]
\[
M_f = K_f \cap Z(\mc S)(\A_f),
\]
\[
\Gam_{K_f} = K_f \cap G(\Q).
\]
If $X$ is the symmetric space for $\mc G(\R)$, then
\begin{equation}\label{eq:cusp-centralizer}
e(\Gam_{K_f} \bs X) \geq c(Z(\mc S), M_f) = c(\mc H, \hat M_f),
\end{equation}
where $\mc H = Z(\mc S) / \mc S$ and $\hat M_f$ is the image of $M_f$ in $\mc H$. Moreover, we have equality in \eqref{eq:cusp-centralizer} when $\mc G$ has an Iwasawa decomposition with respect to $K_f$ and $\mc P$.
\end{prop}


\begin{pf}
We will prove that
\[
c(\mc P, L_f) = c(Z(\mc S), M_f).
\]
Since the torus $\mc S$ is the multiplicative group over $\Q$, it has class number one. Therefore, the right-hand equality in \eqref{eq:cusp-centralizer} follows from Lemma \ref{lem:central-ext}. The proposition then follows immediately from Lemmas \ref{lem:cusp-bound} and \ref{lem:cusp-iwa}.

Let $\mc U$ be the unipotent $\Q$-group such that $\mc P = Z(\mc S) \mc U$. By Corollary 2.5 and Proposition 2.7 in \cite{Borel1}, it suffices to show that
\[
L_f^\infty = M_f^\infty N_f^\infty,
\]
where $N_f^\infty = \mc U(\R) \times (K_f \cap \mc U(\A_f))$, and it suffices to prove the analogous decomposition at any nonarchimedean place. However, when the component of $K_f$ at a fixed nonarchimedean place is a parahoric subgroup, this decomposition follows immediately from the italicized statement in \cite[$\S$3.1.1]{Tits2}. This proves the proposition.
\end{pf}


Lastly, we will need the following relationship between class numbers of unitary and special unitary groups. This result may be known in greater generality, but we could not find a reference.


\begin{prop}\label{prop:unitary-class}
Let $D$ be $\Q$, an imaginary quadratic field, or a definite quaternion algebra over $\Q$, and let $\tau$ be trivial, the nontrivial Galois automorphism, or quaternion conjugation, respectively. Let $h$ be a $\tau$-hermitian form on $D^N$, $\mc H$ the pin/unitary group of $h$, and $\mc H_0 \subset \mc H$ the spin/special unitary group. Let $M_f \subset \mc H(\A_f)$ be an open compact subgroup, $M_f^\infty = \mc H(\R) \times M_f \subset \mc H(\A)$, and $L_f^\infty = \mc H_0(\A) \cap M_f^\infty$. Then there is a universal constant $c$ so that
\begin{equation}\label{eq:unitary-class}
\#\left( M_f^\infty \bs \mc H(\A) / \mc H(\Q) \right) \geq \frac{1}{c} \# \left( L_f^\infty \bs \mc H_0(\A) / \mc H_0(\Q) \right)
\end{equation}
\end{prop}


\begin{pf}
It suffices to show that for any $g \in \mc H_0(\A)$, the elements of $\mc H_0(\A)$ in the double coset $M_f^\infty g \mc H(\Q)$ project to at most $c$ elements of $L_f^\infty \bs \mc H(\A) / \mc H(\Q)$. Suppose that $g_1 = k g h$ for $g_1 \in \mc H_0(\A)$, $k \in M_f^\infty$, and $h \in \mc H(\Q)$. Since $\mc H_0(\A)$ is the kernel of the determinant map
\[
d : \mc H(\A) \to D^*(\A),
\]
we have $d(k) = d(h)^{-1}$. The image of $d$ is contained in the subgroup $D^1(\A)$ of elements in $D(\A)$ of reduced $\tau$-norm $1$. Since $d(M_f^\infty)$ lies in an open compact subgroup of $D^1(\A)$ and $d(\mc H(\Q))$ lies in the rational points, it follows that $d(k)$ and $d(h)$ must lie in $d(M_f^\infty) \cap d(\mc H(\Q))$, which is contained in the subgroup $\mc O^1$ of reduced $\tau$-norm $1$ elements of some $\Z$-order $\mc O$ of $D(\Q)$.

Then $\mc O^1$ is a finite group of order bounded by a universal constant $c_0$. Indeed, $c_0 = 2$ if $D = \Q$, is $6$ is $D$ is imaginary quadratic, and is $24$ when $D$ is a definite quaternion algebra. Therefore, we can choose elements $r_1, \dots, r_{c_0}$ in $M_f^\infty$ and $s_1, \dots, s_{c_0} \in \mc H(\Q)$ so that if $k \in M_f^\infty$ and $d(k) = d(r_j)$ then there exists $\ell \in L_f^\infty$ so that $k = \ell r_j$ and if $d(h) = d(s_j)$ for some $h \in \mc H(\Q)$ then there exists $h_0 \in \mc H_0(\Q)$ so that $h = s_j h_0$.

Therefore, if $g_1 \in \mc H_0(\A)$ and $g_1 = k g h$ for some $k \in M_f^\infty$ and $h \in \mc H(\Q)$, there exist $r_i$, $s_j$, $\ell \in L_f^\infty$, and $h_0 \in \mc H_0(\Q)$ such that
\[
g_1 = \ell (r_i g s_j) h_0.
\]
Then $r_i g s_j \in \mc H_0(\A)$, so $g_1$ has the same image in $L_f^\infty \bs \mc H_0(\A) / \mc H_0(\Q)$ as one of the elements of the finite set $\{r_i g s_j\}$. There are at most $c = c_0^2$ such elements, so this proves the proposition.
\end{pf}


Now we explain these results in greater detail for real and complex hyperbolic space. In particular, we give sufficient information to compute the number of cusps for principal arithmetic subgroups of interest.

\bigskip

\noindent
\textbf{Hyperbolic space}. Recall from $\S$\ref{subsec:real}, $\mc H_0$ is the spin group of $q'$ where $q = q_0 \oplus q'$ for $q_0$ a hyperbolic plane. As explained in \cite{Borel2} for the orthogonal group, the class number of $\mc H_0$ is the number of classes in the spinor genus of $q'$ with respect to the lattice $\mc L'$, where $\mc L'$ is the summand of $\mc L$ associated with $q'$. In particular, we see that Theorem \ref{thm:intro} follows for principal arithmetic lattices from the fact that there are only finitely many anisotropic quadratic forms over $\Q$ with at most $n$ classes in its spinor genus \cite[Appendix A.3]{Cassels}.

For the Bianchi groups $\SL_2(\mc O_k)$, the number of cusps is well-known to equal the class number $h_k$ of the imaginary quadratic field $k$. See \cite{Chinburg--Long--Reid} for a proof. The quadratic form over $\Q$ that determines the corresponding subgroup of $\Spin(3, 1) \cong \SL_2(\C)$ is the binary quadratic form of discriminant equal to the discriminant of the imaginary field $k$. In particular, the above methods show that $\SL_2(\mc O_k)$ has $h_k$ cusps via Gauss's work on binary quadratic forms.

\bigskip

\noindent
\textbf{Complex hyperbolic space}. Let $\ell$ be an imaginary quadratic field and $h$ a hermitian form of signature $(n, 1)$ on the $\ell$-vector space $V$ of dimension $n + 1$. Then $\mc G$ is the special unitary group of $h$ and $\mc H$ is the unitary group of the anisotropic summand $h'$. Therefore, the number of cusps correspond precisely to the class number of the unitary group of a hermitian lattice. This is closely related to the class number and restricted class number of imaginary quadratic fields. See \cite{Stover} for a complete analysis in the case $n = 2$ and \cite{Zeltinger, Kato} for special cases in higher dimensions.

\bigskip

Quaternionic hyperbolic space has a similar description in terms of a hermitian form on a definite quaternion algebra over $\Q$. We leave it to the interested reader to work out this case and the Cayley hyperbolic plane in further detail. We now proceed to maximal lattices and the proof of Theorem \ref{thm:intro}.


\section{Maximal lattices}\label{sec:maximal}


Let $\mc G$ be a simply connected $\Q$-algebraic group of real and $\Q$-rank one as in $\S$\ref{sec:principal}, and let $\Lam < \mc G(\R)$ be a maximal lattice. To prove Theorem \ref{thm:intro}, we must show that for every $x \in \mathbb N$, there are only finitely many $\mc G$ such that $e(\Gam \bs X) \leq x$, where $X$ is the symmetric space associated with $\mc G$. From \cite[Prop.~1.4]{Borel--Prasad}, we know that there exists a coherent open compact subgroup $K_f \subset \mc G(\A_f)$ such that $\Lam$ is the normalizer in $\mc G(\R)$ of $\Gam_{K_f} = \mc G(\Q) \cap K_f$.

Let $\mc S$ be a maximal $\Q$-split torus of $\mc G$, and $\mc H = Z(\mc S) / \mc S$ be the quotient of the centralizer of $\mc S$ in $\mc G$ by $\mc S$. Let $\mc H_0$ be the subgroup of $\mc H$ consisting of elements with determinant one. It follows from Propositions \ref{prop:cusp-centralizer} and \ref{prop:unitary-class} that
\begin{equation}\label{eq:normal-cusp}
e(\Lam \bs X) \geq \frac{e(\Gam_{K_f} \bs X)}{[\Lam : \Gam_{K_f}]} \geq \frac{1}{c} \frac{c(\mc H_0, \hat M_f)}{[\Lam : \Gam_{K_f}]},
\end{equation}
where $\hat M_f$ is the open compact subgroup of $\mc H_0(\A_f)$ determined by $K_f$ and $c$ is the constant from Proposition \ref{prop:unitary-class}. Also, recall from $\S$\ref{sec:principal} that $\mc H_0$ uniquely determines $\mc G$. Therefore, it suffices to show that the right-hand side of \eqref{eq:normal-cusp} is bounded above by $x$ for only finitely many $\mc H_0$.

Let $\mc C$ be the center of $\mc G$. Then \cite[Prop.~2.9]{Borel--Prasad} gives an exact sequence
\[
1 \to \mc C(\R) / (\mc C(\Q) \cap \Gam_{K_f}) \to \Lam / \Gam_{K_f} \to \delta(\overline{\mc G}(\Q))'_\Theta \to 1,
\]
where $\overline{\mc G}$ is the adjoint form of $\mc G$ and $\delta(\overline{\mc G}(\Q))'_\Theta$ is the image of $\Lam$ in $\mathrm H^1(\Q, \mc C)$. The central elements of $\Lam$ clearly act trivially on $X_\infty$, so we in fact have
\begin{equation}\label{eq:galois-cusps}
e(\Lam \bs X) \geq \frac{1}{c} \frac{c(\mc H_0, \hat M_f)}{\# \delta(\overline{\mc G}(\Q))'_\Theta}.
\end{equation}
Once we prove that the right-hand side of \eqref{eq:galois-cusps} is bounded above by $x$ for only finitely many $\mc H_0$, the proof of Theorem \ref{thm:intro} will be complete. In other words, we need to prove the following.


\begin{thm}\label{thm:apply-prasad}
Let $\mc G$ be a $\Q$-algebraic group of real and $\Q$-rank one. Choose a maximal $\Q$-split torus $\mc S$ in $\mc G$, and let $Z(\mc S)$ be the centralizer of $\mc S$ in $\mc G$. Let $\mc H$ be the $\Q$-algebraic group $Z(\mc S) / \mc S$ and $\mc H_0$ be the subgroup of elements with determinant one. Let $\mc C$ be the center of $\mc G$. If $K_f$ is a coherent open compact subgroup of $\mc G(\A_f)$, let $\hat M_f$ be the induced coherent open compact subgroup of $\mc H_0(\A_f)$. For any $x \in \mathbb N$, there are only finitely many such $\mc G$ with
\[
\frac{c(\mc H_0, \hat M_f)}{\# \delta(\overline{\mc G}(\Q))'_\Theta} \leq x,
\]
where $\delta(\overline{\mc G}(\Q))'_\Theta$ denotes the image of $\overline{\mc G}(\Q)$ in $\mathrm H^1(\Q, \mc C)$ and $\overline{\mc G}$ is the adjoint form of $\mc G$.
\end{thm}


\begin{pf}
Applying Prasad's volume formula \cite{Prasad} as in $\S$7.4 of \cite{Borel--Prasad}, we see that
\[
c(\mc H_0, \hat M_f) \geq \mathrm D_\ell^{\frac{1}{2} \mathfrak s(\mc H')} \left( \prod_{i = 1}^r \frac{m_i!}{(2 \pi)^{m_i + 1}} \right) \zeta(\hat M_f).
\]
Here, $\mc H'$ is the unique quasi-split simply connected inner form of $\mc H_0$, $r$ is the absolute rank of $\mc H_0$, $\ell$ is the unique extension of $\Q$ over which $\mc H'$ splits, $\mathrm D_\ell$ is the absolute discriminant of $\ell$, $\mathfrak s(\mc H')$ and the $m_i$ are as in \cite{Prasad} or \cite[$\S$3.7]{Borel--Prasad}, and $\zeta(\hat M_f)$ is the Euler product associated with $\hat M_f$ as in \cite[$\S$7.4]{Borel--Prasad}. Note that $\ell$ is either $\Q$ or a quadratic field and that the Tamagawa number of $\mc H_0$ is one.

Note that $\mc G$ has absolute rank $r + 2$, since $\mc H'$ has absolute rank $r$. Simultaneously considering the three cases arising in Propositions 5.1 and 5.6 of \cite[$\S$5]{Borel--Prasad}, we have the bound
\[
\# \delta(\overline{\mc G}(\Q))'_\Theta \leq 2 h_\ell \mathrm D_\ell n^{2 + \#\mathrm T} \prod_{v < \infty} \# \Xi_{\Theta_v}.
\]
Here $h_\ell$ is the class number of $\ell$, $n$ is as defined in \cite[$\S$2.6]{Borel--Prasad} (so $n \leq 4$ unless $\mathcal G$ is of type $\mathrm A_{r + 2}$, in which case $n = r + 3$), $\mathrm T$ is the set of places of $\ell$ defined before Proposition 5.1 or in $\S$5.5 of \cite{Borel--Prasad}, and $\Xi_{\Theta_v}$ is defined in \cite[$\S$2]{Borel--Prasad}. Note that in combining all the cases, we have introduced overkill into several of them.

This gives
\begin{equation}\label{eq:final-inequal}
\frac{c(\mc H_0, \hat M_f)}{\# \delta(\overline{\mc G}(\Q))'_\Theta} \geq \frac{\mathrm D_\ell^{\frac{1}{2} \mathfrak s(\mc H') - 1}}{2 h_\ell n^{2 + \#\mathrm T}} \left( \prod_{i = 1}^r \frac{m_i!}{(2 \pi)^{m_i + 1}} \right) \prod_{v < \infty} \left( \# \Xi_{\Theta_v} \right)^{-1} \zeta(\hat M_f).
\end{equation}
To prove that the right-hand side of \eqref{eq:final-inequal} is bounded above by $x$ for only finitely many $\mc H_0$ (thus finitely many $\mc G$), it suffices to prove the following two claims:
\begin{eqnarray}\label{eq:last-two}
\frac{1}{n^{\#\mathrm T}} \prod_{v < \infty} \left( \# \Xi_{\Theta_v} \right)^{-1} \zeta(\hat M_f) &>& \del_0 \qquad \textrm{with}\ \del_0\ \textrm{independent of}\ \hat M_f \label{eq:parahoric} \\
\frac{\mathrm D_\ell^{\frac{1}{2} \mathfrak s(\mc H') - 1}}{2 h_\ell n^2} \left( \prod_{i = 1}^r \frac{m_i!}{(2 \pi)^{m_i + 1}} \right) &\geq& x / \delta_0 \qquad \textrm{for finitely many}\ \mc H \label{eq:brauer-siegel}
\end{eqnarray}
Our proofs closely resemble the proof of Theorem B in \cite{Borel--Prasad}, though we emphasize that these inequalities are not an immediate consequence of the computations in \cite{Borel--Prasad}. Indeed, they are considering a single group, whereas each of \eqref{eq:parahoric} and \eqref{eq:brauer-siegel} has factors from both $\mc G$ and $\mc H_0$.

We first prove \eqref{eq:parahoric}. It suffices to show that for all but finitely many $\hat M_f$ and finitely many nonarchimedean places $v$ of $\Q$, we have
\begin{equation}\label{eq:local-euler}
\frac{e(\hat M_v)}{n^{\epsilon(v)} \# \Xi_{\Theta_v}} > 1,
\end{equation}
where $\epsilon(v) \in \{1, 0\}$ depending of whether or not $v \in \mathrm T$, $\hat M_v$ denotes the component of $\hat M_f$ at $v$, and $\zeta(\hat M_f) = \prod_v e(\hat M_v)$ are the Euler factors of $\zeta(\hat M_f)$. When $M_v$ is a special parahoric subgroup and $v \notin \mathrm T$, then $\# \Xi_{\Theta_v} = 1$ and the inequality follows from the fact that $e(\hat M_v) > 1$. For all other $v$, we have the inequality
\begin{equation}\label{eq:euler-bound}
e(\hat M_v) > \frac{p^{r_v + 1}}{p_v + 1},
\end{equation}
where $p_v$ is the rational prime associated with $v$ and $r_v$ is the local rank at $v$. See \cite[$\S$7.4(4)]{Borel--Prasad}.

For the other places, $n^{\epsilon(v)} \# \Xi_{\Theta_v} \leq (r + 3)^2$ (see $\S$2.6 and $\S$A.7 of \cite{Borel--Prasad}). Therefore \eqref{eq:local-euler} holds for all $r_v \geq 8$, independent of $p_v$, and for $p_v \geq 17$, independent of $r$. For each of these finitely many $r$, we compute the product over the finitely many $p_v$ such that $p^{r_v + 1} / (p_v + 1) < 1$ and see that $\del_0 = 0.015$ suffices for all $\mc H_0$. This proves \eqref{eq:parahoric}.

Now consider \eqref{eq:brauer-siegel}. Since $\ell$ is $\Q$ or a quadratic field, the Brauer--Siegel Theorem (more precisely, its proof) implies that
\begin{equation}\label{eq:br-sieg}
h_\ell \leq \left( \frac{5 \pi}{6} \right)^2 \mathrm D_\ell.
\end{equation}
See \cite[$\S$6]{Borel--Prasad}. This implies that
\[
\frac{1}{2 h_\ell n^2} \mathrm D_\ell^{\frac{1}{2} \mathfrak s(\mc H') - 1} \left( \prod_{i = 1}^r \frac{m_i!}{(2 \pi)^{m_i + 1}} \right) \geq
\]
\[
\left( \frac{6}{5 \pi} \right)^2 \frac{\mathrm D_\ell^{\frac{1}{2} \mathfrak s(\mc H') - 2}}{2 n^2} \left( \prod_{i = 1}^r \frac{m_i!}{(2 \pi)^{m_i + 1}} \right).
\]
For \cite[$\S$3.7]{Borel--Prasad}, $\mathfrak s(\mc H') \geq 5$ when $\mc H'$ is not split. When $\mc H'$ is split, $\ell = \Q$ and the discriminant and class number terms are automatically $1$. In particular, independent of $\ell$, we see that
\begin{equation}\label{eq:dell}
\mathrm D_\ell^{\frac{1}{2} \mathfrak s(\mc H') - 2} \geq 1
\end{equation}
for all $\mc H'$. Indeed, when $\mathrm D_\ell$ isn't already $1$, the exponent is at least $0$.

To prove that the right-hand side of \eqref{eq:brauer-siegel} is bounded above by $x / \del_0$ for only finitely many $\mc H'$, it now suffices to show that
\[
\frac{1}{2 n^2} \prod_{i = 1}^r \frac{m_i!}{(2 \pi)^{m_i + 1}}
\]
is bounded above for only finitely many $\mc H'$. This is very similar to the proof of Proposition 6.1 in \cite{Borel--Prasad}. The numbers $\{m_i\}$ form a nondecreasing sequence. Moreover $m_i = m_{i + 1}$ can only occur once within the sequence, and only happens for groups of type $\mathrm D_r$ for $r$ even. Since $m_r \to \infty$ linearly in $r$ \cite[$\S$1.5]{Prasad}, $n \leq r + 3$, and since there are only finitely many $\mc H'$ of bounded absolute rank (since $\mc H'$ is the unique quasi-split form), the result follows from Stirling's Formula.

It remains to prove that finiteness is independent of $\mc H_0$. That is, for a fixed $\mc H'$, we need to know that there are only finitely many $\mc H_0$ with quasi-split simply connected inner form $\mc H'$ such that the class number of $\mc H_0$ with respect to some coherent open compact subgroup of $\mc H_0(\A_f)$ is bounded above. However, this follows from the fact, proven in $\S$7 of \cite{Borel--Prasad}, that $\zeta(\hat M_f)$ is bounded above by any real number $x$ for only finitely many equivalence classes of coherent products $\hat M_f$ of parahoric subgroups (not to mention for only finitely many $\mc H_0$). This completes the proof of the theorem, hence of Theorem \ref{thm:intro}.
\end{pf}


\begin{rem}
Fix $k \in \mathbb N$ and let $\mc G$ and $\mc H_0$ be as above. There are, independent of $\mc H_0$, only finitely many $\mc H_0(\A_f)$-conjugacy classes of $\hat M_f \subset \mc H_0(\A_f)$ such that $c(\mc H_0, \hat M_f)$ is bounded by $k$. This does not imply for a fixed $\mc G$ that there are only finitely many conjugacy classes of $K_f \subset \mc G(\A_f)$ so that $e(\Gam_{K_f} \bs X) \leq k$. Indeed, infinitely many $K_f \subset \mc G(\A_f)$ could induce the same $\hat M_f \subset \mc H_0(\A_f)$. That is, a fixed commensurability class could contain infinitely many distinct \emph{minimal} orbifolds with $\leq k$ ends. This is indeed the case for $\SL_2(\R)$, $\SL_2(\C)$, and $\SU(2, 1)$, where one can build infinite families of distinct but commensurable minimal one-cusped orbifolds. See \cite{Chinburg--Long--Reid, Stover} for the construction.
\end{rem}


\section{One-cusped hyperbolic orbifolds}\label{sec:onecusp}


We now use the techniques of $\S$\ref{sec:maximal} to prove Theorem \ref{thm:intro1cusp}. Let $\Lam$ be a maximal nonuniform arithmetic lattice in $\Spin(n, 1)$ with associated algebraic group $\mc G$. Then $\mc G$ is the spin group of the quadratic form $q$ on $\Q^{n + 1}$ and $\mc H_0$ is the spin group of $q'$, where $q = q_0 \oplus q'$ and $q_0$ is a hyperbolic plane. As above, $\mc H_0$ is the subgroup of elements with determinant one in $\mc H = \mathrm{Pin}(q')$, which is the group $Z(\mc S) / \mc S$, where $\mc S$ is a maximal torus of $\mc G$ and $Z(\mc S)$ its centralizer.

Note that we can take $c = 4$ in Proposition \ref{prop:unitary-class}, since $\mc H_0$ has index two in $\mc H$. Therefore, to show that $e(\Lam \bs \mathbf H^n) > 1$, it suffices by \eqref{eq:galois-cusps} to show that
\begin{equation}\label{eq:hyp-inequal}
\frac{1}{4} \frac{c(\mc H_0, \hat M_f)}{\# \delta(\overline{\mc G}(\Q))'_\Theta} > 1
\end{equation}
for all $n \geq 30$, where the terminology in \eqref{eq:hyp-inequal} is all described in $\S$\ref{sec:maximal}. Since $\mc G$ has type $\mathrm B$ or $\mathrm D$ depending on the parity of $n$, we prove \eqref{eq:hyp-inequal} in two steps.


\begin{pf}[Proof of \eqref{eq:hyp-inequal} for $n$ even]
Here $\mc G$ and $\mc H_0$ have type $\mathrm B_s$ for $s = n / 2, (n - 2) / 2$. We then have the inequality
\[
\# \delta(\overline{\mc G}(\Q))'_\Theta \leq 2^{1 + \# \mathrm T} \prod_v \# \Xi_{\Theta_v},
\]
where $\mathrm T$ is the finite set of nonarchimedean places over which $\mc G$ doesn't split (see \cite[$\S$3.5]{Belolipetsky}). Recall from $\S$\ref{sec:maximal} that we also have
\[
c(\mc H_0, \hat M_f) \geq \left( \prod_{i = 1}^r \frac{m_i!}{(2 \pi)^{m_i + 1}} \right) \zeta(\hat M_f),
\]
where $r = (n - 2) / 2$ is the absolute rank of $\mc H_0$, since $\mc H_0$ is a spin group in $n - 1$ variables, $\zeta(\hat M_f)$ is the Euler factor for $\hat M_f$, and $m_i = 2 i - 1$.

We saw in the proof of Theorem \ref{thm:apply-prasad} that
\[
\frac{1}{2^{\# \textrm T}} \prod_v \left( \# \Xi_{\Theta_v} \right)^{-1} \zeta(\hat M_f) > 1.
\]
as long as every local rank as at least $8$. Therefore, to show that $e(\Lam \bs \mathbf H^n) > 1$, it suffices to show that
\[
\left( \prod_{i = 1}^r \frac{m_i!}{(2 \pi)^{m_i + 1}} \right) > 8.
\]
This product is less than $1$ for $n = 28$ and greater than $9$ for $n = 30$, so this completes the proof of Theorem \ref{thm:intro1cusp} for $n$ even.
\end{pf}


\begin{pf}[Proof of \eqref{eq:hyp-inequal} for $n$ odd]
Now our groups $\mc G$ and $\mc H_0$ have type $\mathrm D_s$ for $s = (n \pm 1) / 2$. The cases $s$ even and $s$ odd are slightly different, so we treat each separately.

Suppose that $\mc H_0$ has type $\mathrm D_r$ for $r$ even. Then $\mc G$ has type ${}^1 \mathrm D_{r + 1}$ and
\[
\# \delta(\overline{\mc G}(\Q))'_\Theta \leq 4^{\# \mathrm T} \prod_v \# \Xi_{\Theta_v}
\]
by \cite[Prop.~4.12(2)(a)]{Belolipetsky--Emery}. Again we have
\[
\frac{1}{2^{\# \textrm T}} \prod_v \left( \# \Xi_{\Theta_v} \right)^{-1} \zeta(\hat M_f) > 1
\]
for rank at least $8$ from the proof of Theorem \ref{thm:apply-prasad}. Therefore, it suffices to show that
\begin{equation}\label{eq:dr-eq-odd}
\left( \prod_{i = 1}^r \frac{m_i!}{(2 \pi)^{m_i + 1}} \right) > 4.
\end{equation}
Here,
\begin{equation}\label{eq:dr-prod}
\left( \prod_{i = 1}^r \frac{m_i!}{(2 \pi)^{m_i + 1}} \right) = \frac{(r - 1)!}{(2 \pi)^r} \left( \prod_{i = 1}^{r - 1} \frac{(2 i - 1)!}{(2 \pi)^{2 i}} \right),
\end{equation}
and a direct check shows that \eqref{eq:dr-eq-odd} holds for all even $r \geq 16$, i.e., for all $n \geq 33$ that are congruent to $1$ modulo $4$.

Now suppose that $\mc H_0$ has type $\mathrm D_r$ for $r$ odd. Then, \cite[Prop.~4.12(2)(c)]{Belolipetsky--Emery}, where we take $\mathrm T = \mathrm R \cup \mathrm T_1$, we have
\[
\# \delta(\overline{\mc G}(\Q))'_\Theta \leq 4^{1 + \# \mathrm T} h_\ell \prod_v \# \Xi_{\Theta_v},
\]
where $\ell$ is a quadratic field. Then
\[
\frac{1}{4^{\# \textrm T}} \prod_v \left( \# \Xi_{\Theta_v} \right)^{-1} \zeta(\hat M_f) > 1
\]
for $r$ at least $8$. Therefore, applying \eqref{eq:br-sieg} and \eqref{eq:dell} we want to know when $r$ is sufficiently large that
\begin{equation}\label{eq:dr-eq-even}
\left( \frac{6}{5 \pi^2} \right) \left( \prod_{i = 1}^r \frac{m_i!}{(2 \pi)^{m_i + 1}} \right) > 8.
\end{equation}
The product on the left-hand side of \eqref{eq:dr-eq-even} is the same as \eqref{eq:dr-prod} (note that $r - 1$ now appears twice), and we see that \eqref{eq:dr-eq-even} holds for all $r \geq 17$, i.e., for any $n \geq 35$ that is congruent to $3$ modulo $4$.

It remains to rule out the case $n = 31$. Here $\mc H_0$ has type $\mathrm D_{15}$. Also factoring
\[
\mathrm D_\ell^{\mf s(\mc H_0) - 2} \geq 3^{\mf s(\mc H_0) - 2} = 3^{\frac{31}{2}}
\]
into \eqref{eq:dr-eq-even} gives the necessary bound and completes the proof of Theorem \ref{thm:intro1cusp}.
\end{pf}


As mentioned in the introduction, there are known one-cusped hyperbolic $n$-orbifolds for all $n \leq 9$ \cite{Hild}. It is not known if there is a one-cusped hyperbolic manifold for any $n \geq 4$. Using the above methods, one could ostensibly build one-cusped orbifolds for $10 \leq n \leq 29$ using definite quadratic forms $q'$ with few classes in their spinor genus. We now build examples of one-ended orbifolds for $n = 10, 11$. We do not know if these appear elsewhere in the literature.

Definite rational quadratic forms with one class in their genus only exist in dimensions below $10$ \cite{Watson}. Furthermore, these are precisely the forms with one class in their spinor genus \cite{Earnest--Hsia}. There are, up to adjoints, $2$ in $9$ variables and $1$ in $10$ variables \cite{Watson910}. If $q'$ is any such form and $q$ is the direct sum of $q'$ and a hyperbolic plane, then $q$ has signature $(n, 1)$ for $n = 10, 11$ and therefore determines a commensurability class of noncompact arithmetic hyperbolic $n$-orbifolds. If $\mc G$ is the $\Q$-algebraic group $\Spin(q)$, choose $K_f \subset \mc G(\A_f)$ so that $\mc G$ has an Iwasawa decomposition with respect to $K_f$ and such that $\mc H_0 = \Spin(q')$ has class number one with respect to the induced open compact subgroup $\hat M_f$ of $\mc H_0$. Then $\mc H = \mathrm{Pin}(q')$ also has class number one with respect to $M_f$, since $\mc H(\Q)$ contains a coset representative for $\mc H_0(\A)$ in $\mc H(\A)$. In particular, Proposition \ref{prop:cusp-centralizer} implies that $e(\Gam_{K_f} \bs \mathbf H^n) = 1$.

It is possible that there are one-cusped arithmetic hyperbolic $n$-orbifolds for $12 \leq n \leq 29$. However, there does not seem to be a classification of the two-class definite rational quadratic forms, so we do not know if Theorem \ref{thm:intro1cusp} is sharp. Given such a form, there is a principal arithmetic lattice $\Gam_{K_f}$ with $2$ cusps by Proposition \ref{prop:cusp-centralizer}, so there is more work needed in order to show that the form determines a one-ended orbifold. One must then show that the normalizer $\Lam$ of $\Gam_{K_f}$ in $\Spin(n, 1)$ is a nontrivial extension and that the covering group $\Lam / \Gam_{K_f}$ for the regular covering $\Gam_{K_f} \bs \mathbf H^n \to \Lam \bs \mathbf H^n$ identifies the two ends of $\Gam_{K_f} \bs \mathbf H^n$.

We close with one final remark regarding one-cusped orbifolds. Instead of using arithmetic orbifolds, one might try to use the construction of Gromov and Piatetski-Shapiro \cite{GPS} to build hybrid examples. In dimension $n \geq 5$, codimension one totally geodesic subspaces of noncompact arithmetic hyperbolic $n$-orbifolds are well-known to be arithmetic, but they are also noncompact. It follows from Theorem \ref{thm:intro1cusp} that one also cannot use a hybrid construction to build one-cusped orbifolds for $n \geq 31$, since the totally geodesic submanifold along which one must cut already has too many cusps, not to mention whatever cusps are on either side of the hypersurface. Similarly, one cannot use this hybrid construction to build a $k$-cusped hyperbolic $n$-orbifold in cases where all arithmetic hyperbolic $(n - 1)$-orbifolds have more than $k$ cusps.


\bibliography{classnumber}


\end{document}